\documentclass[12pt]{article}
\usepackage{amsfonts}
\usepackage{amsmath}
\usepackage{amssymb}
\usepackage{amsthm}
\newtheorem{theorem}{Theorem}
\newtheorem{cor}[theorem]{Corollary}

\newenvironment{proof'}{\textbf{Proof}}{\hfill $\square$ \\}

\newcommand{\inprod}[2]{\langle #1 \, ,#2 \rangle}

\newcommand{\sph}{S^{n-1}}

\title{There are Significantly More Nonnegative Polynomials than Sums of Squares}

\date{September 2003}
\begin{document}
\author{Grigoriy Blekherman}
\numberwithin{theorem}{section} \numberwithin{equation}{theorem}
\maketitle
\begin{abstract}
We investigate the quantitative relationship between nonnegative
polynomials and sums of squares of polynomials. We show that if
the degree is fixed and the number of variables grows then there
are significantly more nonnegative polynomials than sums of
squares. More specifically, we take compact bases of the cone of
nonnegative polynomials and the cone of sums of squares and derive
bounds for the volumes of the bases. If the degree is greater than
2 then we show that the ratio of the volumes of the bases, raised
to the power reciprocal to the ambient dimension, tends to 0 as
the number of variables tends to infinity.
\end{abstract}
\section{Introduction}
\hspace{.44cm} Let $P_{n,2k}$ be the vector space of real
homogeneous polynomials in $n$ variables of degree $2k$. The study
of the relationship between nonnegative polynomials in $P_{n,2k}$
and sums of squares of polynomials of degree $k$ was initiated by
Hilbert. He showed that for the cases $n=2$, $k=1$ and, $n=3$ and
$k=2$, a nonnegative polynomial is necessarily a sum of squares;
in all other cases there exist nonnegative polynomials that are
not sums of squares. Hilbert's proof of existence of nonnegative
polynomials that are not sums of squares was non-constructive
\cite{hilbert}. The first such explicit polynomials were
constructed only fifty years later by Motzkin in the 1940's. There
are currently several known families of
non-negative polynomials that are not sums of squares \cite{choi},\cite{rez1}. \\
\indent There remains however a natural question, which we call,
the \textit{Quantitative Sums of Squares Problem}:
\begin{center}
Are there significantly more nonnegative polynomials \\
than sums of squares of polynomials?
\end{center}

\indent The known examples of nonnegative polynomials that are not
sums of squares live either on the boundary of the cone of
non-negative polynomials or very close to it. The difficulty in
constructing such polynomials explicitly leads naturally to asking
whether they are pathological examples, while the bulk of
nonnegative polynomials
are sums of squares.\\
\indent In this paper we show that if the degree $2k$ is fixed and
we let the number of variables grow then there are significantly
more nonnegative polynomials than sums of squares. One could hope
that sums of squares take up a perhaps small but constant portion
of the cone of nonnegative polynomials, but we will show this is
not the case. As a corollary it follows that if the number of
variables is large then nonnegative polynomials that are not sums
of squares are far more ``normal''
than sums of squares themselves.\\
\indent The Quantitative Sums of Squares Problem also has
significance from the point of view of computational complexity.
Using the tools of semidefinite programming one can efficiently
compute whether a given polynomial is a sum of squares
\cite{parrilo}. However, determining whether a polynomial is
nonnegative is NP-hard for $k \geq 2$ \cite[Part 1]{complexity}.
Therefore we can ask how much we lose by testing for sums of
squares instead of nonnegativity. There has been some experimental
evidence that for certain subsets of the cone of nonnegative
polynomials sums of squares approximate nonnegative polynomials
rather well \cite{parrilo}. However, it follows from the results
of this paper that if the number of variables is large compared to
the degree then there are far fewer sums of squares than
nonnegative polynomials, and thus relaxation of testing for
nonnegativity to testing for sums of squares
does not work well.\\

\section{Main Theorems}
\hspace{4.4mm} We begin by introducing some notation. Nonnegative
polynomials and sums of squares form full-dimensional convex cones
in $P_{n,2k}$. Let $C(=C_{n,2k})$ be the cone of nonnegative
polynomials,
\begin{equation*}
C=\bigl{\{}f \in P_{n,2k} \mid f(x) \geq 0 \quad \text{for all}
\quad x \in \mathbb{R}^n \bigr{\}}.
\end{equation*}
Let $Sq(=Sq_{n,2k})$ be the cone of sums of squares,
\begin{equation*}
Sq=\biggl{\{} f \in P_{n,2k} \mathrel{\bigg{\arrowvert}} f=\sum_i
f_i^2 \quad \text{for some} \quad f_i \in P_{n,k} \biggr{\}}.
\end{equation*}
We work with the following inner product on $P_{n,2k}$,
\begin{equation*}
\inprod{f}{g}=\int_{S^{n-1}} fg \, d\sigma,
\end{equation*}
where $\sigma$ is the rotation invariant probability measure on
$\sph$.\\
\indent In order to compare the cones $C$ and $Sq$ we take compact
bases. Let $M(=M_{n,2k})$ be the hyperplane of all forms in
$P_{n,2k}$ with integral 0 on the unit sphere $\sph$:
\begin{equation*}
M=\biggr{\{} f \in P_{n,2k} \mathrel{\bigg{\arrowvert}}
\int_{\sph}f \, d\sigma=0 \biggl{\}}.
\end{equation*}
We use $D_M$ to denote the dimension of $M$,  $S_M$ to denote the
unit sphere in $M$ and $B_M$ to denote the unit ball in $M$.\\
\indent Let $r^{2k}$ in $P_{n,2k}$ be the following polynomial
\begin{equation*}
r^{2k}=(x_1^2+ \ldots +x_n^2)^k.
\end{equation*}
We define compact convex bodies $\widetilde{C}$ and
$\widetilde{Sq}$ as the sets of all forms $f$ in $M$ such that
$f+r^{2k}$ lies in the respective cone:
\begin{equation*}
\widetilde{C}=\{f \in M \mid f+r^{2k} \in C \}, \qquad
\widetilde{Sq}=\{f \in M \mid f+r^{2k} \in Sq \},
\end{equation*}
Convex bodies $\widetilde{C}$ and $\widetilde{Sq}$ are sections of
their respective cones with the hyperplane of forms of integral
one, which are translated into $M$ by subtracting $r^{2k}$.\\
\indent The main result of this paper are the two estimates below:
\begin{theorem}
\label{infinitynormbound} There is the following lower bound on
the volume of $\widetilde{C}$:
\begin{equation*}
\bigg{(}\frac{\text{Vol} \hspace{.8mm} \widetilde{C}}{\text{Vol}
\, B_M}\bigg{)}^{\frac{1}{D_M}} \, \geq \frac{1}{2\sqrt{4k+2}} \,
n^{-1/2}.
\end{equation*}
\end{theorem}
\begin{theorem}
\label{squarenormbound} There is the the following upper bound on
the volume of $\widetilde{Sq}$:
\begin{equation*}
\bigg{(}\frac{\text{Vol} \, \widetilde{Sq}}{\text{Vol} \,
B_M}\bigg{)}^{\frac{1}{D_M}} \leq
\frac{4^{2k}(2k)!\sqrt{24}}{k!}\, n^{-k/2}.
\end{equation*}
\end{theorem}
From the above two theorems we trivially obtain the following
corollary which allows us to compare $C$ and $Sq$:
\begin{cor}
\label{maincor} There exists a constant $c(k)$ dependent only on
$k$ such that
\begin{equation*}
\bigg{(}\frac{\text{Vol} \hspace{.8mm} \widetilde{C}}{\text{Vol}
\, \widetilde{Sq}}\bigg{)}^{\frac{1}{D_M}} \, \geq c(k)
n^{(k-1)/2}.
\end{equation*}
It suffices to take
\begin{equation*}
c(k)=\frac{k! }{2(2k)!4^{2k}\sqrt{24(4k+2)}}.
\end{equation*}
\end{cor}
We note that if $k=1$ then $\widetilde{C}=\widetilde{Sq}$ and
indeed the bound of Corollary \ref{maincor} becomes just $c(1)$.
However, if we fix $k \geq 2$ then there are clearly far more
nonnegative polynomials than sums of squares.

\section{Volumes of Nonnegative Polynomials}
\hspace{.44cm} In this section we prove Theorem
\ref{infinitynormbound}. For a real vector space $V$ with the unit
sphere $S_V$ and a function $f:V \rightarrow \mathbb{R}$ we use
$||f||_{p}$ to denote the $L^p$ norm of $f$:
\begin{equation*}
||f||_p=\bigg{(}\int_{S_V} |f|^p \, d\mu \bigg{)}^{1/p} \qquad
\text{and} \qquad ||f||_{\infty}=\max_{x \in S_V} |f(x)|.
\end{equation*}
\indent We begin by observing that $\widetilde{C}$ is a convex
body in $M_{n,2k}$ with origin in its interior and the boundary of
$\widetilde{C}$ consists of polynomials with minimum $-1$ on
$\sph$. Therefore the gauge $G_C$ of $\widetilde{C}$ is given by:
\begin{equation*}
G_{C}(f)=|\min_{v \, \in \sph} f(v)\,|.
\end{equation*}
By using integration in polar coordinates in $M$ it is easy to
obtain the following expression for the volume of $\widetilde{C}$,
\begin{equation}
\label{integralvolume} \biggl{(}\frac{\text{Vol} \,
\widetilde{C}}{\text{Vol }
B_M}\biggr{)}^{\frac{1}{D_M}}=\biggl{(}\int_{S_M} G_{C}^{-D_M} \,
d\mu \biggr{)}^{\frac{1}{D_M}},
\end{equation}
where $\mu$ is the rotation invariant probability measure on
$S_M$. The relationship \eqref{integralvolume}
holds for any convex body with origin in its interior \cite[p. 91]{pisier}. \\
We interpret the right hand side of \eqref{integralvolume} as
$||G_C^{-1}||_{D_M}$, and by H\"{o}lder's inequality
\begin{equation*}
||G_C^{-1}||_{D_M} \geq ||G_C^{-1}||_1.
\end{equation*}
Thus,
\begin{equation*}
\biggl{(}\frac{\text{Vol} \, \widetilde{C}}{\text{Vol} \,
B_M}\biggr{)}^{\frac{1}{D_M}} \geq \int_{S_M} G_C^{-1} \, d\mu.
\end{equation*}
By applying Jensen's inequality \cite[p.150]{hardy}, with convex
function $y=1/x$ it follows that,
\begin{equation*}
\int_{S_M} G_C^{-1} \ d\mu \geq \bigg{(}\int_{S_M} G_C \, d\mu
\bigg{)}^{-1}.
\end{equation*}
Hence we see that
\begin{equation*}
\bigg{(} \frac{\text{Vol} \, \widetilde{C}}{\text{Vol} \, B_M}
\bigg{)}^{\frac{1}{D_M}} \geq \bigg{(}\int_{S_M} |\min f| \, d \mu
\bigg{)}^{-1}.
\end{equation*}
Clearly, for all $f \in P_{n,2k}$
\begin{equation*}
||f||_{\infty} \geq |\min f|.
\end{equation*}
Therefore,
\begin{equation*}
\bigg{(} \frac{\text{Vol} \, \widetilde{C}}{\text{Vol} \, B_M}
\bigg{)}^{\frac{1}{D_M}} \geq \bigg{(}\int_{S_M} ||f||_{\infty} \,
d \mu \bigg{)}^{-1}.
\end{equation*}
The proof of Theorem \ref{infinitynormbound} is now completed by
the following estimate.
\begin{theorem}
\label{infinitynorm} Let $S_M$ be the unit sphere in $M$ and let
$\mu$ be the rotation invariant probability measure on $S_M$. Then
the following inequality for the average $L^{\infty}$ norm over
$S_M$ holds:
\begin{equation*}
\int_{S_M} ||f||_{\infty} \, d\mu \leq 2\sqrt{2n(2k+1)}.
\end{equation*}
\end{theorem}
\begin{proof}
It was shown by Barvinok in \cite{barv} that for all $f \in
P_{n,2k}$,
\begin{equation*}
||f||_{\infty} \leq \binom{2kn+n-1}{2kn}^{\frac{1}{2n}}||f||_{2n}.
\end{equation*}
By applying Stirling's formula we can easily obtain the bound
\begin{equation*}
\binom{2kn+n-1}{2kn}^{\frac{1}{2n}} \leq 2\sqrt{2k+1}.
\end{equation*}
Therefore it suffices to estimate the average $L^{2n}$ norm, which
we denote by $A$:
\begin{equation*}
A=\int_{S_M} ||f||_{2n} \, d\mu.
\end{equation*}
Applying H\"{o}lder's inequality we observe that
\begin{equation*}
A=\int_{S_M} \bigg{(}\int_{\sph} f^{2n}(x) \, d\sigma
\bigg{)}^{\frac{1}{2n}} d\mu \leq \bigg{(}\int_{S_M}\int_{\sph}
f^{2n}(x)\, d\sigma \, d\mu \bigg{)}^{\frac{1}{2n}}.
\end{equation*}
By interchanging the order of integration we obtain
\begin{equation}
\label{first} A \leq \bigg{(} \int_{\sph} \int_{S_M} f^{2n}(x) \,
d\mu \, d\sigma \bigg{)}^{\frac{1}{2n}}.
\end{equation}
We now note that by symmetry of $M$
\begin{equation*}
\int_{S_M} f^{2n}(x) \, d\mu,
\end{equation*}
is the same for all $x \in \sph$. Therefore we see that in
\eqref{first} the outer integral is redundant and thus
\begin{equation}
\label{second} A \leq \bigg{(} \int_{S_M} f^{2n}(v) \, d\mu
\bigg{)}^{\frac{1}{2n}}, \qquad \text{where} \ v \ \text{is any
vector in} \ \sph.
\end{equation}
\indent For $v \in \mathbb{R}^n$ the functional
\begin{equation*}
\lambda_v:M \longrightarrow \mathbb{R}, \qquad \lambda_v(f)=f(v),
\end{equation*}
is linear and therefore there exists a form $q_v \in M$ such that
\begin{equation*}
\lambda_v(f)=\inprod{q_v}{f}.
\end{equation*}
Rewriting \eqref{second} we see that
\begin{equation}
\label{third} A \leq \bigg{(} \int_{S_M} \inprod{f}{q_v}^{2n} \,
d\mu \bigg{)}^{\frac{1}{2n}}.
\end{equation}
There are explicit descriptions of the polynomials $q_v$, see for
example \cite{vilenkin}, we will only need the property that for
$v \in \sph$
\begin{equation*}
||\,q_v||_{2}=\sqrt{D_M}.
\end{equation*}
This can also be deduced by abstract representation theoretic
considerations.
\\\indent We observe that
\begin{equation*}
\int_{S_M} \inprod{f}{q_v}^{2n} \, d\mu=(D_M)^n \
\frac{\Gamma(n+\frac{1}{2}\,) \,
\Gamma(\frac{1}{2}D_M)}{\sqrt{\pi} \, \Gamma(\frac{1}{2}D_M+n)}.
\end{equation*}
We substitute this into \eqref{third} to obtain,
\begin{equation*}
A \leq \bigg{(}(D_M)^n \ \frac{\Gamma(n+\frac{1}{2}\,) \,
\Gamma(\frac{1}{2}D_M)}{\sqrt{\pi} \, \Gamma(\frac{1}{2}D_M+n)}
\bigg{)}^{\frac{1}{2n}}.
\end{equation*}
Since
\begin{equation*}
\bigg{(}\frac{\Gamma(\frac{1}{2}D_M)}{\Gamma(\,\frac{ 1}{2}D_M
+n)}\bigg{)}^{\frac{1}{2n}} \leq \sqrt{\frac{2}{D_M}} \qquad
\text{and} \qquad \bigg{(}
\frac{\Gamma(n+1/2\,)}{\sqrt{\pi}}\bigg{)}^{\frac{1}{2n}} \leq
n^{1/2},
\end{equation*}
we see that
\begin{equation*}
A \leq (2n)^{1/2}.
\end{equation*}
The theorem now follows.
\end{proof}
\section{Volumes of Sums of Squares}
 \setcounter{equation}{0}
In this section we prove Theorem \ref{squarenormbound}. Let us
begin by considering the support function of $\widetilde{Sq}$,
which we call $L_{\widetilde{Sq}}$:
\begin{equation*}
L_{\widetilde{Sq}}(f)=\max_{g \, \in \, \widetilde{Sq}} \,
\inprod{f}{g}.
\end{equation*}
The average width $W_{\widetilde{Sq}}$ of $\widetilde{Sq}$ is
given by
\begin{equation*}
W_{\widetilde{Sq}}=2\int_{S_M} L_{\widetilde{Sq}} \, d\mu.
\end{equation*}
We now recall Urysohn's Inequality \cite[p.318]{schneider} which
applied to $\widetilde{Sq}$ gives
\begin{equation}
\label{urysohn} \bigg{(}\frac{\text{Vol} \,
\widetilde{Sq}}{\text{Vol} \, B_M}\bigg{)}^{\frac{1}{D_M}} \, \leq
\frac{W_{\widetilde{Sq}}}{2}.
\end{equation}
Therefore it suffices to obtain an upper bound for
$W_{\widetilde{Sq}}$. \\
\indent Let $S_{P_{n,k}}$ denote the unit sphere in $P_{n,k}$. We
observe that extreme points of $\widetilde{Sq}$ have the form
\begin{equation*}
g^2-r^{2k} \qquad \text{where} \qquad g \in P_{n,k} \qquad
\text{and} \qquad \int_{\sph}g^2 \,d\sigma=1.
\end{equation*}
For $f \in M$,
\begin{equation*}
\inprod{f}{r^{2k}}=\int_{\sph}f \,d\sigma=0,
\end{equation*}
and therefore,
\begin{equation*}
L_{\widetilde{Sq}}(f)=\max_{g \, \in S_{P_{n,k}}} \inprod{f}{g^2}.
\end{equation*}
\indent We now introduce a norm on $P_{n,2k}$, which we denote $||
\ ||_{sq}$:
\begin{equation*}
||f||_{sq}=\max_{g \, \in \, S_{P_{n,k}}} |\inprod{f}{g^2}|.
\end{equation*}
It is clear that
\begin{equation*}
L_{Sq}(f) \leq ||f||_{Sq}.
\end{equation*}
Therefore by \eqref{urysohn} it follows that
\begin{equation*}
\bigg{(}\frac{\text{Vol} \, \widetilde{Sq}}{\text{Vol} \,
B_M}\bigg{)}^{\frac{1}{D_M}} \, \leq \int_{S_M} ||f||_{sq} \, d
\mu.
\end{equation*}
The proof of Theorem \ref{squarenormbound} is reduced to the
estimate below.
\begin{theorem}
There is the following bound for the average $|| \ ||_{sq}$ over
$S_M$:
\begin{equation*}
\int_{S_M} ||f||_{sq} \, d\mu \, \leq \,
\frac{4^{2k}(2k)!\sqrt{24}}{k!}\, n^{-k/2}.
\end{equation*}
\end{theorem}
\begin{proof}
For $f \in P_{n,2k}$ we introduce a quadratic form $H_f$ on
$P_{n,k}$:
\begin{equation*}
H_f(g)=\inprod{f}{g^2} \qquad \text{for} \qquad g \in P_{n,k}.
\end{equation*}
We note that
\begin{equation*}
||f||_{sq}=\max_{g \, \in \,
S_{P_{n,k}}}|\inprod{f}{g}|=||H_f||_{\infty}.
\end{equation*}
We bound $||H_f||_{\infty}$ by a high $L^{2p}$ norm of $H_f$.
Since $H_f$ is a form of degree 2 on the vector space $P_{n,k}$ of
dimension $D_{n,k}$ it follows by the inequality of Barvinok in
\cite{barv} applied in the same way as in the proof of Theorem
\ref{infinitynormbound} that
\begin{equation*}
||H_f||_{\infty} \leq 2\sqrt{3} \, ||H_f||_{2D_{n,k}}.
\end{equation*}
Therefore it suffices to estimate:
\begin{equation*}
A= \int_{S_M} ||H_f||_{2D_{n,k}} \, d\mu = \int_{S_M}
\bigg{(}\int_{S_{P_{n,k}}} \inprod{f}{g^2}^{\, 2D_{n,k}} \, d
\sigma(g) \, d\mu(f) \bigg{)}^{\frac{1}{2D_{n,k}}}.
\end{equation*}
We apply H\"{o}lder's inequality to see that
\begin{equation*}
A \leq \bigg{(}\int_{S_M} \int_{S_{P_{n,k}}} \inprod{f}{g^2}^{\,
2D_{n,k}} \, d \sigma(g) \, d\mu(f) \bigg{)}^{\frac{1}{2D_{n,k}}}.
\end{equation*}
By interchanging the order of integration we obtain
\begin{equation}
\label{weird1} A \leq \bigg{(}\int_{S_{P_{n,k}}} \int_{S_M}
\inprod{f}{g^2}^{\, 2D_{n,k}} \, d \mu(f) \, d\sigma(g)
\bigg{)}^{\frac{1}{2D_{n,k}}}.
\end{equation}
\indent Now we observe that the inner integral
\begin{equation*}
\int_{S_M} \inprod{f}{g^2}^{\, 2D_{n,k}} \, d \mu(f),
\end{equation*}
clearly depends only on the length of the projection of $g^2$ into
$M$. Therefore we have
\begin{equation*}
\int_{S_M} \inprod{f}{g^2}^{\, 2D_{n,k}} \, d \mu(f) \leq \,
||g^2||_2^{2D_{n,k}}\int_{S_M}\inprod{f}{p}^{2D_{n,k}} \, d\mu(f)
\quad \text{for any} \quad p \in S_M.
\end{equation*}
We observe that
\begin{equation*}
||g^2||_2=(||g||_4)^2 \qquad \text{and} \qquad ||g||_2=1.
\end{equation*}
By a result of Duoandikoetxea \cite{duo} Corollary 3 it follows
that
\begin{equation*}
||g^2||_2 \leq 4^{2k}.
\end{equation*}
Hence we obtain
\begin{equation*}
\int_{S_M} \inprod{f}{g^2}^{\, 2D_{n,k}} \, d \mu(f) \leq
4^{4kD_{n,k}} \int_{S_V}\inprod{f}{p}^{2D_{n,k}} \, d\mu(f).
\end{equation*}
We note that this bound is independent of $g$ and substituting
into \eqref{weird1} we get
\begin{equation*}
A \leq 4^{2k}\bigg{(}\int_{S_V}\inprod{f}{p}^{2D_{n,k}} \,
d\mu(f)\bigg{)}^{\frac{1}{2D_{n,k}}}.
\end{equation*}
\indent Since $p \in S_M$ we have
\begin{equation*}
\int_{S_M}\inprod{f}{p}^{2D_{n,k}} \, d\mu(f) =
\frac{\Gamma(D_{n,k}+\frac{1}{2})\Gamma(\,
\frac{1}{2}D_M)}{\sqrt{\pi} \, \Gamma(D_{n,k}+\frac{1}{2}D_M)}.
\end{equation*}
We use the following easy inequalities:
\begin{equation*}
\bigg{(}\frac{\Gamma(\,
\frac{1}{2}D_M)}{\Gamma(D_{n,k}+\frac{1}{2}D_M)}\bigg{)}^{\frac{1}{2D_{n,k}}}
 \leq \sqrt{\frac{2}{D_M}}
\end{equation*}
and
\begin{equation*}
\bigg{(}\frac{\Gamma(D_{n,k}+\frac{1}{2})}{\sqrt{\pi}}\bigg{)}^{\frac{1}{2D_{n,k}}}
\leq \sqrt{D_{n,k}},
\end{equation*}
to see that
\begin{equation*}
A \leq 4^{2k}\sqrt{\frac{2D_{n,k}}{D_M}}.
\end{equation*}
We now recall that
\begin{equation*}
D_{n,k}=\binom{n+k-1}{k} \qquad \text{and} \qquad
D_M=\binom{n+2k-1}{2k}-1.
\end{equation*}
Therefore
\begin{equation*}
\sqrt{\frac{D_{n,k}}{D_M}} \, \leq  \, \frac{(2k)!}{k!} \,
n^{-k/2}.
\end{equation*}
Thus
\begin{equation*}
A \leq \frac{4^k(2k)!\sqrt{2}}{k!} \, n^{-k/2}.
\end{equation*}
The theorem now follows.
\end{proof}

\section{Remarks}
The estimates of Theorems \ref{infinitynormbound} and
\ref{squarenormbound} for the volumes of nonnegative polynomials
and sums of squares are asymptotically exact if the degree is
fixed. The proof, however, is more technical and shall be
discussed elsewhere.

\medskip
\textsc{Department of Mathematics, University of Michigan, Ann
Arbor, MI 48109-1109, USA}\\
\textit{Email address:} gblekher@umich.edu
\end{document}